\pdfoutput=1
\RequirePackage{ifpdf}
\ifpdf 
\documentclass[pdftex]{sigma}
\else
\documentclass{sigma}
\fi

\begin{document}


\renewcommand{\PaperNumber}{065}

\FirstPageHeading

\ShortArticleName{Special Functions of Hypercomplex Variable on the Lattice Based on SU(1,1)}

\ArticleName{Special Functions of Hypercomplex Variable\\
on the Lattice Based on SU(1,1)}

\Author{Nelson FAUSTINO}

\AuthorNameForHeading{N.~Faustino}

\Address{Departamento de Matem\'atica Aplicada, IMECC--Unicamp,
\\
CEP 13083--859, Campinas, SP, Brasil}
\Email{\href{mailto:faustino@ime.unicamp.br}{faustino@ime.unicamp.br}}
\URLaddress{\url{https://sites.google.com/site/nelsonfaustinopt/}}

\ArticleDates{Received May 06, 2013, in f\/inal form October 28, 2013; Published online November 05, 2013}

\Abstract{Based on the representation of a~set of canonical operators on the lattice~$h\mathbb{Z}^n$, which
are Clif\/ford-vector-valued, we will introduce new families of special functions of hypercomplex variable
possessing $\mathfrak{su}(1,1)$ symmetries.
The Fourier decomposition of the space of Clif\/ford-vector-valued polynomials with respect to the ${\rm
SO}(n)\times \mathfrak{su}(1,1)$-module gives rise to the construction of new families of polynomial
sequences as eigenfunctions of a~coupled system involving forward/backward discretizations $E_h^{\pm}$ of
the Euler operator $E=\sum\limits_{j=1}^nx_j \partial_{x_j}$.
Moreover, the interpretation of the one-parameter representation $\mathbb{E}_h(t)=\exp(tE_h^--tE_h^+)$ of
the Lie group ${\rm SU}(1,1)$ as a~semigroup $\left(\mathbb{E}_h(t)\right)_{t\geq 0}$ will allows us to
describe the polynomial solutions of an homogeneous Cauchy problem on $[0,\infty)\times h{\mathbb Z}^n$
involving the dif\/ferencial-dif\/ference operator $\partial_t+E_h^+-E_h^-$.}

\Keywords{Clif\/ford algebras; f\/inite dif\/ference operators; Lie algebras}

\Classification{22E70; 30G35; 33C80; 39A12}

\section{Introduction}

In the investigation of special functions, the representation theory through Lie groups allows to compute
families of orthogonal polynomials in terms of hypergeometric series expansions (see,
e.g.,~\cite{VilenkinKlimyk91}).
Recent approaches towards discrete quantum mechanics, as for instance in~\cite{OR11}, reveals that the
representation of f\/inite dif\/ference operators as canonical generators of a~certain Lie algebra provides
a~general scheme to construct sequences of polynomials as eigenfunctions of a~discrete Hamiltonian operator.
These sequences of polynomials that appear on the literature under the name of Shef\/fer sequences
(cf.~\cite{DHS96,LTW04}) or Appell sequences (cf.~\cite{Tempesta08}) give rise to families of Bernoulli and
Euler polynomials beyond the classical rising/falling factorial polynomials (see also~\cite[Section~3]{Cartier00}).

Seen the fact that Clif\/ford algebras of signature $(0,n)$ encode the structure of the special orthogonal
group ${\rm SO}(n)$ of $n\times n$ matrices (cf.~\cite[Subsection~I.1]{DSS92}), it remains natural to study
multi-variable extensions of the above approaches to the Euclidean space ${\mathbb R}^n$ in terms of
hypercomplex variables.
In~\cite[Section~2]{MT08} and in~\cite[Section~3]{RSKS10} the authors obtained the
hypercomplex extension of Bernoulli and Euler polynomials, respectively; In~\cite[Section~4]{RSS12}
the authors considered discrete versions of Fueter polynomials as an alternative hypercomplex extension for
the raising/lowering Clif\/ford-vector-valued polynomials considered in~\cite[Section~3]{FK07}.
In~\cite[Section~2]{FR11} the authors shown that such families of Clif\/ford-vector-valued
polynomials of discrete variable may be realized from Lie algebraic representations of an algebra of
endomorphisms analogue to the radial algebra representation obtained in {\it continuum} by Sommen
(cf.~\cite{Sommen97}).

This paper is organized as follows: In Section~\ref{ScopeProblemsSection} it will be given
the motivation to study, in the framework of Clif\/ford algebras, sequences of polynomials generated from
a~set of f\/inite dif\/ference operators.
In Section~\ref{PolynomialSU11} the construction of irreducible representations for the spaces of
Clif\/ford-vector-valued polynomials on the lattice $h{\mathbb Z}^n$ based on the Howe dual pair $({\rm
SO}(n),\mathfrak{su}(1,1))$ (cf.~\cite{Howe89TransAMS}) will be considered.

To f\/ind the Fourier decomposition for the space of Clif\/ford-vector-valued polynomials on the lattice
$h{\mathbb Z}^n$, we will start to determine the positive and negative series representations for the Lie
group ${\rm SU}(1,1)$ in interconnection with the forward/backward discretizations $E_h^\pm$ of the
classical Euler operator $E=\sum\limits_{j=1}^nx_j \partial_{x_j}$.
Afterwards, the action of the ${\rm SO}(n)\times \mathfrak{su}(1,1)$-module on the lattice $h{\mathbb Z}^n$
will produce a~sequence of invariant and irreducible subspaces under the discrete series representations of
${\rm SU}(1,1)$.

The results obtained in Section~\ref{PolynomialSU11} will be used in
Section~\ref{FamiliesSpecialFunctions} to describe the space of Clif\/ford-vector-valued polynomials
on the lattice $h{\mathbb Z}^n$ as hypergeometric series expansions and to characterize the homogeneous
solutions of the dif\/ferential-dif\/ference operator $\partial_t+E_h^+-E_h^-$ in terms of the semigroup
$(\mathbb{E}_h(t))_{t\geq 0}$ carrying the one-parameter representation
$\mathbb{E}_h(t)=\exp(tE_h^--tE_h^+)$ of ${\rm SU}(1,1)$.

\section{Scope of problems}
\label{ScopeProblemsSection}

Let ${\bf e}_1,{\bf e}_2,\ldots,{\bf e}_n$ be an orthogonal basis of ${\mathbb R}^n$.
The Clif\/ford algebra of signature $(0,n)$, which we will denote by $C \kern -0.1em \ell_{0,n}$,
corresponds to the algebra generated from the set of graded anti-commuting relations
\begin{gather}
\label{CliffordBasis}
{\bf e}_j{\bf e}_k+{\bf e}_k{\bf e}_j=-2\delta_{j,k}
\qquad
\text{for any}
\quad
j,k=1,2,\ldots,n.
\end{gather}
Under the linear space isomorphism given by the mapping ${\bf e}_{j_1}{\bf e}_{j_2}\cdots {\bf e}_{j_r}
\mapsto dx_{j_1}dx_{j_2}\cdots dx_{j_r}$, with $1\leq j_1<j_2<\dots<j_r\leq n$, the resulting algebra with
dimension $2^n$ is isomorphic to the exterior algebra $\bigwedge ({\mathbb R}^n)$.

This allows us to represent any vector $x=(x_1,x_2,\ldots,x_n)$ of ${\mathbb R}^n$ as an element
$x=\sum\limits_{j=1}^n x_j {\bf e}_j\in C \kern -0.1em \ell_{0,n}$ and the translations $(x_1,x_2,\ldots,
x_j\pm h,\ldots,x_n)$ on the grid $h{\mathbb Z}^n \subset {\mathbb R}^n$ with mesh width $h>0$ by the
displacements $x\pm h{\bf e}_j$ over $C \kern -0.1em \ell_{0,n}$.
Moreover, the Clif\/ford-vector-valued functions correspond to linear combinations in terms of the
$r$-multivector basis ${\bf e}_{j_1}{\bf e}_{j_2}\cdots {\bf e}_{j_r}$ labeled by subsets
$J=\{j_1,j_2,\ldots,j_r\}$ of $\{ 1,2,\ldots,n\}$, i.e.\
\begin{gather*}
{\bf f}(x)=\sum_{r=0}^n\sum_{|J|=r}f_J(x){\bf e}_J
\qquad
\text{with}
\quad
{\bf e}_{J}={\bf e}_{j_1}{\bf e}_{j_2}\cdots{\bf e}_{j_r}.
\end{gather*}

Let us now recall some basic facts about f\/inite dif\/ference operators.
The forward/backward f\/inite dif\/ferences $\partial_h^{\pm j}$ def\/ined on the grid $h{\mathbb Z}^n$ by
\begin{gather*}
\big(\partial_{h}^{+j}{\bf f}\big)(x)=\frac{{\bf f}(x+h{\bf e}_j)-{\bf f}(x)}{h}
\qquad
\text{and}
\qquad
\big(\partial_{h}^{-j}{\bf f}\big)(x)=\frac{{\bf f}(x)-{\bf f}(x-h{\bf e}_j)}{h}
\end{gather*}
are interrelated by the translation operators $\big(T_h^{\pm j}{\bf f}\big)(x)={\bf f}(x\pm h {\bf e}_j)$, i.e.
\begin{gather}
\label{TranlationsPmj}
T_h^{-j}\big(\partial_h^{+j}{\bf f}\big)(x)=\big(\partial_h^{-j}{\bf f}\big)(x)
\qquad
\text{and}
\qquad
T_h^{+j}\big(\partial_h^{-j}{\bf f}\big)(x)=\big(\partial_h^{+j}{\bf f}\big)(x).
\end{gather}
Also, they satisfy the product rules
\begin{gather}
\partial_h^{+j} ({\bf g}(x){\bf f}(x) )=\big(\partial_{h}^{+j}{\bf g}\big)(x){\bf f}(x+h{\bf e}_j)+{\bf g}
(x)\big(\partial_h^{+j}{\bf f}\big)(x),
\nonumber
\\
\partial_h^{-j} ({\bf g}(x){\bf f}(x) )=\big(\partial_{h}^{-j}{\bf g}\big)(x){\bf f}(x-h{\bf e}_j)+{\bf g}
(x)\big(\partial_h^{-j}{\bf f}\big)(x).
\label{productRule}
\end{gather}

Along the paper we will use the bold letters ${\bf f},{\bf g},\ldots,{\bf m},\ldots,{\bf w},\ldots$ and so
on, when we refer to functions of the above form.
The f\/inite dif\/ference Dirac operators $D_h^\pm$ def\/ined {\it viz}
\begin{gather*}
D_h^+=\sum_{j=0}^n{\bf e}_j\partial_h^{+j}
\qquad
\text{and}
\qquad
D_h^-=\sum_{j=0}^n{\bf e}_j\partial_h^{-j}
\end{gather*}
are Clif\/ford-vector-valued and correspond to f\/inite dif\/ference approximations of the classical
gradient operator in a~coordinate-free way.

Now let $\mathcal{P}={\mathbb R}[x]\otimes C \kern -0.1em \ell_{0,n}$ be the space of
Clif\/ford-vector-valued polynomials.
We say that $\{ {\bf m}_s(x;\tau): \, s\in {\mathbb N}_0\}\subset \mathcal{P}$ is an Appell set carrying
$D_h^+$, resp.~$D_h^-$, if for $\tau=\pm h$ we have $D_h^+ {\bf m}_0(x,-h)=D_h^- {\bf m}_0(x;h)=0$ and
\begin{gather}
D_h^+{\bf m}_s(x;-h)=s{\bf m}_{s-1}(x;-h),
\nonumber
\\
D_h^-{\bf m}_s(x;h)=s{\bf m}_{s-1}(x;h),
\qquad
\text{for every}
\ \
s\in{\mathbb N}.\label{AppellSet}
\end{gather}

Iterating $r$-times the operator $D_h^{+}$, resp.\ $D_h^{-}$, it turns out that the action of the semigroup
$\left( \exp(tD_h^+)\right)_{t\geq 0}$, resp.\ $\left( \exp(tD_h^-)\right)_{t\geq 0}$, on each ${\bf
m}_s(x;-h)$, resp.~${\bf m}_s(x;h)$, gives rise to a~binomial expansion.
Indeed, the iterated relations $\left(D_h^{\pm}\right)^r {\bf m}_s(x;\mp h)=\frac{s!}{(s-r)!}{\bf
m}_{s-r}(x;\mp h)$, resp.\ $\left(D_h^{\pm}\right)^r {\bf m}_s(x;\mp h)=0$, that hold for $s\geq r$, resp.\
$s<r$, leads to
\begin{gather*}
\exp\left(tD_h^{+}\right){\bf m}_{s}(x;-h)=\sum_{r=0}^\infty\frac{t^r}{r!}\left(D_h^{+}\right)^r{\bf m}_{s}
(x;-h)=\sum_{r=0}^s\left(
\begin{matrix} s
\\
r
\end{matrix}
\right)t^r{\bf m}_{s-r}(x;-h),
\end{gather*}
and analogously, to $\exp\left(tD_h^{-}\right){\bf m}_{s}(x;h)=\sum\limits_{r=0}^s \left(
\begin{matrix} s
\\
r
\end{matrix}
\right) t^r {\bf m}_{s-r}(x;h)$.

The action of the semigroups $\left( \exp(tD_h^+)\right)_{t\geq 0}$ and $\left( \exp(tD_h^-)\right)_{t\geq
0}$ on $\mathcal{P}$ then correspond to the hypercomplex extension of the Taylor series expansion for
polynomials (see, e.g.,~\cite[Subsection~3.3]{Cartier00}) that gives rise e.g.~to Clif\/ford-vector-valued
polynomials of Bernoulli type analogous to the ones obtained in~\cite[Section~2]{MT08}.

This approach corresponds to the discrete counterpart of the Cauchy--Kovaleskaya extension described
in~\cite[Subsection~III.2]{DSS92}.
For an alternative application of this approach in interrelationship with discrete versions of Fueter
polynomials, one refer to~\cite[Sections~3,~4]{RSS12}.

The Fock space formalism carrying Hilbert spaces (cf.~\cite{BLR98}) reveals that the problem of
constructing polynomial sets $\{ {\bf m}_s(x;\tau):\, s\in {\mathbb N}_0\}$, with $\tau=\pm h$, possessing
the Appell set property is equivalent to the construction of operator endomorphisms $M_h^{+},M_h^-\in
\text{End}(\mathcal{P})$ in such way that the elements of the form
\begin{gather*}
{\bf m}_s(x;\pm h)=\lambda_s\big(M_h^{\pm}\big)^{s}{\bf m}_0(x;\pm h),
\qquad
s\in{\mathbb N}_0
\end{gather*}
yield a~basis for $\mathcal{P}$.
Hereby, the constants $\lambda_s$ are chosen under the condition $\lambda_0=1$ and the
constraints~\eqref{AppellSet}.

In terms of the umbral calculus formalism (see~\cite{DHS96,LTW04} and the references therein), from the
identity operator $I:{\bf f}(x)\mapsto {\bf f}(x)$ and the commuting bracket $[A,B]$ def\/ined as
\begin{gather*}
[A,B]{\bf f}(x)=A(B{\bf f}(x))-B(A{\bf f}(x)),
\end{gather*}
one can start to construct the set of Clif\/ford-vector-valued polynomials $\{ {\bf m}_s(x;\tau): s\in
{\mathbb N}_0\}$ from Weyl--Heisenberg algebra symmetries.
In order to proceed, we will def\/ine for a~given linear polynomial $w(t) \in {\mathbb R}[t]$ of degree $1$
satisfying $(\partial_{h}^{+j}w)(x_k)=(\partial_{h}^{-j}w)(x_k)=\delta_{jk}\mu $, the following set of
multiplication operators:
\begin{gather}
(W_j{\bf f})(x) = \mu^{-1}w(x_j){\bf f}(x),\nonumber
\\
\big(W_h^{+j}{\bf f}\big)(x) = \mu^{-1}w\left(x_j+\frac{h}{2}\right){\bf f}(x+h{\bf e}_j),\nonumber
\\
\big(W_h^{-j}{\bf f}\big)(x) = \mu^{-1}w\left(x_j-\frac{h}{2}\right){\bf f}(x-h{\bf e}_j).\label{LadderWj}
\end{gather}

It is now straightforward from the product rules~\eqref{productRule} that the set of operators
\begin{gather*}
\big\{W_h^{-j},\partial_h^{+j},I : j=1,2,\ldots,n\big\}
\qquad
\text{and}
\qquad
\big\{W_h^{+j},\partial_h^{-j}
,I : j=1,2,\ldots,n\big\}
\end{gather*}
span the Weyl--Heisenberg algebra of dimension $2n+1$.
The remainder graded commuting relations are given by
\begin{alignat}{4}
& \big [\partial_{h}^{+j},\partial_{h}^{+k}\big]=0, \qquad && \big[W_h^{-j},W_h^{-k}\big]=0, \qquad
&& \big[\partial_{h}^{+j},W_h^{-k}\big]=\delta_{jk}I, & \label{WeylHeisenbergWjMinus} \\
\label{WeylHeisenbergWjPlus}
& \big[\partial_{h}^{-j},\partial_{h}^{-k}\big]=0, \qquad
&&
\big[W_h^{+j},W_h^{+k}\big]=0,\qquad
&&
\big[\partial_{h}^{-j},W_h^{+k}\big]=\delta_{jk}I. &
\end{alignat}

In the discrete Clif\/ford analysis setting (cf.~\cite[Subsection~1.3]{FR11}), it is precisely the
Weyl--Heisenberg relations~\eqref{WeylHeisenbergWjMinus}, resp.~\eqref{WeylHeisenbergWjPlus}, that allows us
to determine in a~unique way $M_h^-$, resp.~$M_h^+$, as
\begin{gather*}
M_h^-=\sum_{j=1}^n{\bf e}_j W_h^{-j},
\qquad
\text{resp.}
\qquad
M_h^+=\sum_{j=1}^n{\bf e}_j W_h^{+j},
\end{gather*}
providing in this way a~{\it Fourier duality} between $D_h^+$ and $M_h^-$, resp.~$D_h^-$ and $M_h^+$.

For discretizations of the Dirac operator $D=\sum\limits_{j=1}^n{\bf e}_j \partial_{x_j}$, written as
a~superposition of forward and backward dif\/ferences, the {\it Fourier duality} may be constructed by
means of a~set of skew-Weyl relations (cf.~\cite{RSKS10}).

Notice that the {\it Fourier duality} terminology comes from invariant theory (cf.~\cite{Howe89TransAMS}).
In the language of Clif\/ford analysis this is nothing else than the so-called {\it Fischer duality}
(cf.~\cite{Sommen97}).

Since $M_h^{\pm}$ maps ${\bf m}_s(x;\pm h)$ into ${\bf m}_{s+1}(x;\pm h)$, it remains clear that each
Clif\/ford-vector-valued polynomial ${\bf m}_s(x;h)$, resp.\ ${\bf m}_s(x;-h)$, is an eigenfunction for the
factorized Hamiltonian $M_h^+D_h^-+ D_h^-M_h^+$, resp.\ $M_h^-D_h^++ D_h^+M_h^-$.
The Weyl--Heisenberg character between the operators~$\partial_h^{+j}$ and $W_h^{-j}=\mu^{-1}
w\left(x_j-\frac{h}{2}\right) T_{h}^{-j}$, resp.\ $\partial_h^{-j}$ and $W_h^{+j}=\mu^{-1}
w\left(x_j+\frac{h}{2}\right) T_{h}^{+j}$, combined with the orthogonality constraint~\eqref{CliffordBasis}
provided by the Clif\/ford basis of $C \kern -0.1em \ell_{0,n}$ allows us to rewrite the factorized
Hamiltonians $M_h^+D_h^-+ D_h^-M_h^+$ and $M_h^-D_h^++ D_h^+M_h^-$ in terms of the following
forward/backward discretizations of the Euler operator $E=\sum\limits_{j=1}x_j\partial_{x_j}$:
\begin{gather}
\label{Ehpm}
E_h^{+}=\sum_{j=1}^n\mu^{-1}w\left(x_j+\frac{h}{2}\right)\partial_h^{+j},
\qquad
E_h^{-}=\sum_{j=1}^n\mu^{-1}w\left(x_j-\frac{h}{2}\right)\partial_h^{-j}.
\end{gather}

Indeed, from~\eqref{TranlationsPmj} one can rewrite $E_h^+$ and $E_h^-$ as
\begin{gather*}
E_h^{+}=\sum_{j=1}^n W_h^{+j} \partial_h^{-j}
\qquad
\text{and}
\qquad
E_h^{-}=\sum_{j=1}^n W_h^{-j} \partial_h^{+j}.
\end{gather*}

Then $M_h^+D_h^-+ D_h^-M_h^+$ and $M_h^-D_h^++ D_h^+M_h^-$ admit the following decompositions:
\begin{gather*}
M_h^+D_h^-+D_h^-M_h^+=\sum_{j=1}^n\big({-}2W_h^{+j}\partial_h^{-j}-I\big)
=-2E_h^+-nI,\\
M_h^-D_h^++D_h^+M_h^-=\sum_{j=1}^n\big({-}2W_h^{-j}\partial_h^{+j}-I\big)
=-2E_h^--nI.
\end{gather*}

The second quantization approach through the set of Weyl--Heisenberg
generators~\eqref{WeylHeisenbergWjMinus} and~\eqref{WeylHeisenbergWjPlus} (cf.~\cite{BLR98}) reveals that $E_h^\pm$ are number-type operators with spectrum~${\mathbb N}_0$.
This leads to
\begin{gather*}
E_h^{+}{\bf m}_s(x;h)=s{\bf m}_s(x;h)
\qquad
\text{and}
\qquad
E_h^{-}{\bf m}_s(x;-h)=s{\bf m}_s(x;-h).
\end{gather*}

From the above construction, it turns out that the Weyl--Heisenberg algebra character provided from the
graded commuting relations~\eqref{WeylHeisenbergWjMinus}, resp.\ \eqref{WeylHeisenbergWjPlus}, gives the
required ladder structure to construct also discrete analogues for Gegenbauer polynomials as a~series
expansion written in terms of the Appell set $\{ {\bf m}_s(x;\tau) : s\in {\mathbb N}_0\}$
(cf.~\cite{Eelbode12}).

Since the eigenfunctions of $E_h^+$ and $E_h^-$ do not coincide, in general, it remains natural to ask
which polynomial subsets $\{{\bf m}_s(x;h) : s\in {\mathbb N}_0\}$ of $\mathcal{P}$ give rise to solutions
of the coupled eigenvalue system
\begin{gather*}
E_h^+{\bf m}_s(x;h)=s{\bf m}_s(x;h),
\qquad
E_h^-{\bf m}_s(x;h)=s{\bf m}_s(x;h),
\qquad
s\in{\mathbb N}_0.
\end{gather*}

From the relations~\eqref{TranlationsPmj}, it remains clear that forward and backward dif\/ferences,
$\partial_h^{+j}$ and $\partial_h^{-j}$ respectively, commute.
In contrast, the operators $W_h^{-j}$, $W_h^{+j}$ do not commute in general.
This means that for each $j=1,2,\ldots,n$ the set of operators
$\partial_h^{+j}$, $\partial_h^{-j}$, $W_h^{-j}$, $W_h^{+j}$ and~$I$ do not endow a~canonical realization of an
Weyl--Heisenberg type algebra.

\looseness=1
Although the solutions of the above coupled system of eigenvalue equations may not be represented in terms
of Weyl--Heisenberg algebra symmetries, the set of generators~$W_h^{-j}$, $W_h^{+j}$,~$W_j$ itself (see the
coordinate expressions~\eqref{LadderWj}) will be the departure point of this paper to describe the {hidden}
Lie algebraic symmetries encoded by the solutions of the above coupled system.

\section[Clifford-vector-valued polynomials related with ${\rm SU}(1,1)$]{Clif\/ford-vector-valued
polynomials related with $\boldsymbol{{\rm SU}(1,1)}$}
\label{PolynomialSU11}

\subsection[Discrete series representions of ${\rm SU}(1,1)$]{Discrete series representions of
$\boldsymbol{{\rm SU}(1,1)}$}

Accordingly to~\cite[Section 6.4]{VilenkinKlimyk91}, the Lie group ${\rm SU}(1,1)$ has two families
of discrete series representations.
In order to determine it algebraically, let us take a~close look for the graded commuting relations
involving the number operators $E_h^+$ and $E_h^-$ def\/ined in~\eqref{Ehpm}.

For this purpose, one starts to show that the set of operators $W_h^+$, $W_h^-$ and $W$ def\/ined from the
left endomorphisms~\eqref{LadderWj} acting on the space $\mathcal{P}$:
\begin{gather}
\label{Whpm}
W_h^+=\sum_{j=1}^nW_h^{+j},
\qquad
W_h^-=\sum_{j=1}^nW_h^{-j}
\qquad
\text{and}
\qquad W=\sum_{j=1}^n W_j
\end{gather}
generate a~Lie algebra isomorphic to $\mathfrak{su}(1,1)$.
In order to proceed, we will start with the following lemma which interrelates the set of generators
\begin{gather*}
W_h^{+j}=\mu^{-1}w\left(x_j+\frac{h}{2}\right)T_h^{+j},
\qquad
W_h^{-j}=\mu^{-1}w\left(x_j-\frac{h}{2}\right)T_h^{-j},
\qquad
W_j=\mu^{-1}w(x_j)I.
\end{gather*}
\begin{lemma}
\label{wjCommutingLemma}
For every $j,k=1,2,\ldots,n$ we have the following set of graded commuting rules:
\begin{enumerate}\itemsep=0pt
\item[$(1)$] The operators $W_h^{+j}$ and $W_h^{-k}$ satisfy $\big[W_h^{+j},W_h^{-k}\big]=2h\delta_{jk} W_k.$
\item[$(2)$] The operators $W_h^{+k}$, resp.\ $W_h^{-k}$, and $W_j$ are interrelated by
\begin{gather*}
\big[W_h^{+k},W_j\big]=h\delta_{jk} W_h^{+k},\qquad
\text{resp.}\qquad
\big[W_j,W_h^{-k}\big]=h\delta_{jk} W_h^{-k}.
\end{gather*}
\end{enumerate}
\end{lemma}
\begin{proof}   (1)~From the conditions $T_{h}^{+k} w\left(x_k-\frac{h}{2}\right)=
w\left(x_k+\frac{h}{2}\right)$, $T_{h}^{-k} w\left(x_k+\frac{h}{2}\right)= w\left(x_k-\frac{h}{2}\right)$
and $w\left(x_k \mp \frac{h}{2}\right)=T_{h}^{\pm j}w\left(x_k \mp \frac{h}{2}\right)$, for $j \neq k$, one
obtain for every $j,k=1,2,\ldots,n$, the set of graded commuting relations
\begin{gather*}
\left[w\left(x_j+\frac{h}{2}\right)T_h^{+j},w\left(x_k-\frac{h}{2}\right)T_h^{-k}\right]=\delta_{jk}
\left(w\left(x_k+\frac{h}{2}\right)^2-w\left(x_k-\frac{h}{2}\right)^2\right)I.
\end{gather*}

Now recall that $w(t) \in {\mathbb R}[t]$ is a~polynomial of degree $1$.
Combination of linearity arguments with the condition $\partial_h^{\pm k}w(x_k)=\mu$ lead to the set of
equations $ w\left(x_k+\frac{h}{2}\right)+ w\left(x_k-\frac{h}{2}\right)=2w(x_k)$ and $
w\left(x_k+\frac{h}{2}\right)- w\left(x_k-\frac{h}{2}\right)=h \mu$, and hence, to the set of equations
\begin{gather*}
w\left(x_k+\frac{h}{2}\right)^2-w\left(x_k-\frac{h}{2}\right)^2=2\mu h w(x_k)
\qquad
\text{with}
\quad
k=1,2,\ldots,n.
\end{gather*}
Thus, for all $j,k=1,2,\ldots,n$ the above set of graded commuting relations are equivalent to
$\big[W_h^{+j},W_h^{-k}\big]=2\delta_{jk}h W_k$.

\medskip

(2)~Since $W_k=\mu^{-1}w(x_k)I$ commutes with
\begin{gather*}
W_h^{+j}=\mu^{-1}
w\left(x_j+\frac{h}{2}\right) T_h^{+j}, \qquad \text{resp.}\qquad W_h^{-j}=\mu^{-1} w\left(x_j-\frac{h}{2}\right) T_h^{-j},
\end{gather*}
for every $j \neq k$, it remains to show for every $j=1,2,\ldots,n$ the graded commuting relations
\begin{gather*}
\big[W_h^{+j},W_j\big]=h W_h^{+j}
\qquad
\text{and}
\qquad
\big[W_j,W_h^{-j}\big]=h W_h^{-j}.
\end{gather*}

From a~direct computation
\begin{gather*}
\left[w\left(x_j+\frac{h}{2}\right)T_h^{+j},w(x_j)I\right]=w\left(x_j+\frac{h}{2}\right)^2T_h^{+j}
-w(x_j)w\left(x_j+\frac{h}{2}\right)T_h^{+j}
\\
\phantom{\left[w\left(x_j+\frac{h}{2}\right)T_h^{+j},w(x_j)I\right]}
=h\big(\partial_{h}^{+j}w_j\big)(x) w\left(x_j+\frac{h}{2}\right)T_{h}^{+j},
\\
\left[w(x_j)I,w\left(x_j-\frac{h}{2}\right)T_h^{-j}\right]=w(x_j)w\left(x_j-\frac{h}{2}\right)T_h^{-j}
-w\left(x_j-\frac{h}{2}\right)^2T_h^{-j}
\\
\phantom{\left[w(x_j)I,w\left(x_j-\frac{h}{2}\right)T_h^{-j}\right]}
=h\big(\partial_{h}^{-j}w_j\big)(x) w\left(x_j-\frac{h}{2}\right)T_{h}^{-j}.
\end{gather*}

Combination of the conditions $\big(\partial_{h}^{+j}w\big)(x_j)=\big(\partial_{h}^{-j}w\big)(x_j)=\mu$ with the coordinate
expressions~\eqref{LadderWj} yields $\big[W_h^{+j},W_j\big]=hW_h^{+j}$, resp.\ $\big[W_j,W_h^{-j}\big]=hW_h^{-j}$, as
desired.
\end{proof}

From Lemma~\ref{wjCommutingLemma}, one obtain for the set of multiplication operators~\eqref{LadderWj}, the
Lie algebra isomorphism
\begin{gather*}
\text{span}\left\{\frac{1}{h}W_h^{+j},\frac{1}{h}W_h^{-j},\frac{1}{h}W_j : j=1,2,\ldots,n\right\}
\cong\mathfrak{sl}(2n,{\mathbb R}).
\end{gather*}

Thus, one can infer that $\frac{1}{h}W_h^+$, $\frac{1}{h}W_h^-$ and $\frac{1}{h}W$ are the canonical
generators of the three-dimensional Lie algebra $\mathfrak{su}(1,1)\cong\mathfrak{sl}(2,{\mathbb R})$.
The remaining commuting relations are given by
\begin{gather}
\label{su11Whpm}
\left[\frac{1}{h}W_h^+,\frac{1}{h}W\right]=\frac{1}{h}W_h^+,
\qquad
\left[\frac{1}{h}W_h^-,\frac{1}{h}W\right]=-\frac{1}{h}W,
\qquad
\left[\frac{1}{h}W_h^+,\frac{1}{h}W_h^-\right]=\frac{2}{h}W.
\end{gather}

Now let's turn our attention to the coordinate expressions of $E_h^\pm$ given by~\eqref{Ehpm}.
Recall that the conditions $\big(\partial_{h}^{+ j}w\big)(x_k)=\big(\partial_{h}^{-j}w\big)(x_k)=\delta_{jk}\mu$ provided
from construction allows us to recast~$E_h^+$ and~$E_h^-$ in terms of the multiplication operators
def\/ined in~\eqref{LadderWj}.
Indeed, the basic identities
\begin{gather*}
w\left(x_j+\frac{h}{2}\right)\partial_h^{+j} = \frac{1}{h}\left(w\left(x_j+\frac{h}{2}\right)T_h^{+j}
-w(x_j)I\right)-\frac{\mu}{2}I,
\\
w\left(x_j-\frac{h}{2}\right)\partial_h^{-j} = \frac{1}{h}\left(w(x_j)I-w\left(x_j-\frac{h}{2}
\right)T_h^{-j}\right)-\frac{\mu}{2}I
\end{gather*}
lead to the following coordinate expressions:
\begin{gather}
E_h^{+}=\sum_{j=1}^n\left(\frac{1}{h}W_h^{+j}-\frac{1}{h}W_j-\frac{1}{2}I\right)=\frac{1}{h}
W_h^+-\frac{1}{h}W-\frac{n}{2}I,\nonumber
\\
E_h^{-} = \sum_{j=1}^n\left(\frac{1}{h}W_j-\frac{1}{h}W_h^{-j}-\frac{1}{2}I\right)=\frac{1}{h}W-\frac{1}{h}
W_h^--\frac{n}{2}I,\label{EhpmLadderWj}
\end{gather}
and consequently, to the following coordinate expressions involving the sum/dif\/ference between~$E_h^{+}$
and~$E_h^-$:
\begin{gather}
\label{EhpmLinearCombinations}
E_h^{+}+E_h^{-} = \frac{1}{h}W_h^+-\frac{1}{h}W_h^--n I,
\qquad
E_h^{+}-E_h^{-} = \frac{1}{h}W_h^++\frac{1}{h}W_h^--\frac{2}{h}W.
\end{gather}

The following set of results are also straightforward and will give the key ingredients to construct the
positive/negative series representations encoded by SU(1,1).

\begin{lemma}
\label{su11EhpmLemma}
The operators $E_h^+-E_h^-$, $\frac{1}{h}W_h^+$ and $E_h^++\frac{n}{2}I$, resp.\ $E_h^+-E_h^-$,
$\frac{1}{h}W_h^-$ and $E_h^-+\frac{n}{2}I$ are the canonical generators of the Lie algebra
$\mathfrak{su}(1,1)$.
The remainder commuting relations are given by
\begin{gather*}
\left[E_h^{\pm}+\frac{n}{2}I,E_h^+-E_h^-\right]=E_h^--E_h^+,
\qquad
\left[E_h^{\pm}+\frac{n}{2}I,\frac{1}{h}W_h^\pm\right]=\frac{1}{h}W_h^\pm,
\\
\left[E_h^+-E_h^-,\frac{1}{h}W_h^\pm\right]=2\left(E_h^\pm+\frac{n}{2}I\right).
\end{gather*}
\end{lemma}

\begin{proof} First, notice that direct combination of relations~\eqref{su11Whpm} with the coordinate
expressions~\eqref{EhpmLadderWj} and~\eqref{EhpmLinearCombinations} results into the following set of
graded commuting relations carrying the operators $E_h^+ +\frac{n}{2}I$, $E_h^-+\frac{n}{2}I$ and
$E_h^+-E_h^-$:
\begin{gather*}
\left[E_h^++\frac{n}{2}I,E_h^-+\frac{n}{2}I\right]=\left[\frac{1}{h}W_h^+-\frac{1}{h}W,\frac{1}{h}W-\frac{1}
{h}W_h^-\right]=E_h^+-E_h^-.
\end{gather*}
In order to prove the graded commuting relations $
\left[E_h^{+}{+}\frac{n}{2}I,\frac{1}{h}W_h^+\right]=\frac{1}{h}W_h^+$ and
$\left[E_h^{-}{+}\frac{n}{2}I,\frac{1}{h}W_h^-\right]$ $=\frac{1}{h}W_h^-$ one starts to rewrite
$\left[E_h^++\frac{n}{2}I,\frac{1}{h}W_h^+\right]$ and $\left[ E_h^-+\frac{n}{2}I,\frac{1}{h}W_h^-\right]$
and based on the coordinate expressions~\eqref{EhpmLadderWj}.
In concrete
\begin{gather}
\left[E_h^++\frac{n}{2}I,\frac{1}{h}W_h^+\right] = \left[\frac{1}{h}W_h^+-\frac{1}{h}W,\frac{1}{h}
W_h^+\right],
\nonumber\\
\left[E_h^-+\frac{n}{2}I,\frac{1}{h}W_h^-\right] = \left[\frac{1}{h}W-\frac{1}{h}W_h^-,\frac{1}{h}
W_h^-\right].\label{EhpmCommutator2}
\end{gather}

The graded commuting relations provided from~\eqref{su11Whpm} yield
\begin{gather*}
\left[\frac{1}{h}W_h^+-\frac{1}{h}W,\frac{1}{h}W_h^+\right]=\frac{1}{h}W_h^+
\qquad
\text{and}
\qquad
\left[\frac{1}{h}
W-\frac{1}{h}W_h^-,\frac{1}{h}W_h^-\right]=\frac{1}{h}W_h^-,
\end{gather*}
and therefore, the relations~\eqref{EhpmCommutator2} are equivalent to
\begin{gather*}
\left[E_h^++\frac{n}{2}I,\frac{1}{h}W_h^+\right]=\frac{1}{h}W_h^+
\qquad
\text{and}
\qquad
\left[E_h^-+\frac{n}{2}
I,\frac{1}{h}W_h^-\right]=\frac{1}{h}W_h^-.
\end{gather*}

Finally, the relations $\left[E_h^+\!-E_h^-,\frac{1}{h}W_h^+\right]=2\left( E_h^+\!+\frac{n}{2}I\right)$ and
$\left[E_h^+-E_h^-,\frac{1}{h}W_h^-\right]=2\left( E_h^-+\frac{n}{2}I\right)$ follow straightforwardly from
direct combination of the coordinate expression obtained in~\eqref{EhpmLinearCombinations} for
$E_h^+-E_h^-$ with the graded commutators
\begin{gather*}
\left[\frac{1}{h}W_h^++\frac{1}{h}W_h^--\frac{2}{h}W,\frac{1}{h}W_h^+\right]=2\left(\frac{1}{h}
W_h^+-\frac{1}{h}W\right),
\\
\left[\frac{1}{h}W_h^++\frac{1}{h}W_h^--\frac{2}{h}W,\frac{1}{h}W_h^-\right]=2\left(\frac{1}{h}W-\frac{1}{h}
W_h^-\right).\tag*{\qed}
\end{gather*}
\renewcommand{\qed}{}
\end{proof}

In the proof of Proposition~\ref{su11EhpmProposition} and in the subsequent results, we will make use of
the following lemma which follows straightforwardly from induction on ${\mathbb N}$.
\begin{lemma}
\label{ABsLemma}
For every $A$, $B$ and for every $s\in {\mathbb N}$, the graded commutator $[A,B^s]$ satisfies the summation
formula
\begin{gather*}
[A,B^s]=\sum_{r=0}^{s-1}B^{r}[A,B]B^{s-1-r}.
\end{gather*}
\end{lemma}

\begin{proposition}
\label{su11EhpmProposition}
For any $s \in {\mathbb N}$ we have the following graded commuting relations:
\begin{gather*}
\left[E_h^{\pm}+\frac{n}{2}I,\left(E_h^+-E_h^-\right)^s\right]=-s\left(E_h^+-E_h^-\right)^s,
\\
\left[E_h^{\pm}+\frac{n}{2}I,\left(\frac{1}{h}W_h^{\pm}\right)^s\right]=s\left(\frac{1}{h}W_h^\pm\right)^s,
\\
\left[E_h^+-E_h^-,\left(\frac{1}{h}W_h^\pm\right)^s\right]=s\left(2E_h^\pm+(n-s+1)I\right)\left(\frac{1}{h}
W_h^\pm\right)^{s-1}.
\end{gather*}
\end{proposition}

\begin{proof} Recall that when $[A,B]=\pm B$, Lemma~\ref{ABsLemma} reduces to
\begin{gather*}
\left[A,B^s\right]=\sum_{r=0}^{s-1}\pm B^{r}B B^{s-1-r}=\pm s B^s.
\end{gather*}

Combination of Lemma~\ref{su11EhpmLemma} with the above identity carrying the substitutions
$A=E_h^{\pm}+\frac{n}{2}I$ and $B= E_h^+-E_h^-/B=\frac{1}{h}W_h^\pm$ yield
\begin{gather*}
\left[E_h^{\pm}+\frac{n}{2}I,\left(E_h^+-E_h^-\right)^s\right]=-s\left(E_h^+-E_h^-\right)^s,
\\
\left[E_h^{\pm}+\frac{n}{2}I,\left(\frac{1}{h}W_h^\pm\right)^s\right]=s\left(\frac{1}{h}W_h^\pm\right)^s.
\end{gather*}

For the proof of
\begin{gather*}
\left[E_h^+-E_h^-,\left(\frac{1}{h}W_h^\pm\right)^s\right]=s\left(2E_h^\pm+(n-s+1)I\right)\left(\frac{1}{h}
W_h^\pm\right)^{s-1}
\end{gather*}
recall that the relations $[E_h^\pm +\frac{n}{2}I,\frac{1}{h}W_h^\pm]=\frac{1}{h}W_h^\pm$ provided from
Lemma~\ref{su11EhpmLemma} are equivalent to the intertwining properties $\frac{1}{h}W_h^\pm \left(E_h^\pm
+\frac{n}{2}I\right)=\left(E_h^\pm+\left(\frac{n}{2}-1\right)I\right)\frac{1}{h}W_h^\pm.$

Induction over $r=1,\ldots,s-1$ gives
\begin{gather*}
\left(\frac{1}{h}W_h^\pm\right)^r\left(E_h^\pm+\frac{n}{2}I\right)=\left(E_h^\pm+\left(\frac{n}{2}
-r\right)I\right)\left(\frac{1}{h}W_h^\pm\right)^r.
\end{gather*}

Finally, the direct application of Lemma~\ref{ABsLemma} carrying the substitutions $A=E_h^+-E_h^-$ and
$B=\frac{1}{h}W_h^+$ results into
\begin{gather*}
\left[E_h^+-E_h^-,\left(\frac{1}{h}W_h^\pm\right)^s\right]=\sum_{r=0}^{s-1}\left(\frac{1}{h}
W_h^\pm\right)^r 2\left(E_h^\pm+\frac{n}{2}I\right) \left(\frac{1}{h}W_h^\pm\right)^{s-1-r}
\\
\qquad{} =\sum_{r=0}^{s-1}2\left(E_h^\pm+\left(\frac{n}{2}-r\right)I\right)\left(\frac{1}{h}W_h^\pm\right)^{s-1}
=s\left(2E_h^\pm+(n-s+1)I\right)\left(\frac{1}{h}W_h^\pm\right)^{s-1}.\!\!\!\tag*{\qed}
\end{gather*}
\renewcommand{\qed}{}
\end{proof}

Now let ${\bf m}_0(x;h)$ be a~Clif\/ford-vector-valued polynomial that satisf\/ies the set of equations
$E_h^+{\bf m}_0(x;h)=E_h^-{\bf m}_0(x;h)=0$.
The intertwining relations
$E_h^{\pm}\left(\frac{1}{h}W_h^{\pm}\right)^s=\left(\frac{1}{h}W_h^{\pm}\right)^s(E_h^{\pm}+s I)$ that
yield from
$\left[E_h^{\pm}+\frac{n}{2}I,\left(\frac{1}{h}W_h^\pm\right)^s\right]=s\left(\frac{1}{h}W_h^\pm\right)^s$
(see Proposition~\ref{su11EhpmProposition}) lead to
\begin{gather*}
E_h^+\left[\left(\frac{1}{h}W_h^+\right)^s{\bf m}_0(x;h)\right]=s\left(\frac{1}{h}W_h^+\right)^s{\bf m}
_0(x;h),
\\
E_h^-\left[\left(\frac{1}{h}W_h^-\right)^s{\bf m}_0(x;h)\right]=s\left(\frac{1}{h}W_h^-\right)^s{\bf m}
_0(x;h).
\end{gather*}
Then it is straightforward to see that the basis functions of the form ${\bf
w}_{s}(x;h)=\left(\frac{1}{h}W_h^+\right)^s {\bf m}_0(x;h)$ and ${\bf
w}_{s}(x;-h)=\left(\frac{1}{h}W_h^-\right)^s {\bf m}_0(x;h)$ satisfy the following set of ladder operator
relations:
\begin{gather}
\frac{1}{h}W_h^+ {\bf w}_{s}(x;h)={\bf w}_{s+1}(x;h),\nonumber
\\
\big(E_h^+-E_h^-\big){\bf w}_{s}(x;h)=s(s+n+1){\bf w}_{s-1}(x;h),\label{su11positive}
\\
\left(E_h^++\frac{n}{2}I\right){\bf w}_{s}(x;h)=\left(s+\frac{n}{2}\right){\bf w}_{s}(x;h);\nonumber
\\
\frac{1}{h}W_h^- {\bf w}_{s}(x;-h)={\bf w}_{s+1}(x;-h),\nonumber
\\
\big(E_h^+-E_h^-\big){\bf w}_{s}(x;-h)=s(s+n+1){\bf w}_{s-1}(x;-h),\label{su11negative}
\\
\left(E_h^-+\frac{n}{2}I\right){\bf w}_{s}(x;-h)=\left(s+\frac{n}{2}\right){\bf w}_{s}(x;-h).\nonumber
\end{gather}

We are now in conditions to construct the positive, resp.\ negative, part for the discrete series
representation carrying the group ${\rm SU}(1,1)$ in the same order of ideas of~\cite[Section~6.4]{VilenkinKlimyk91}.
Recall that for the given set of generators $\frac{1}{h}W_h^+$, $\frac{1}{h}W_h^-$ and $\frac{1}{h}W$ of
$\mathfrak{su}(1,1)$, the graded commuting relations~\eqref{su11Whpm} endow a~$*$-structure def\/ined viz
$\left(\frac{1}{h}W\right)^*=\frac{1}{h}W$ and $\left(\frac{1}{h}W_h^+\right)^*=-\frac{1}{h}W_h^-$.
A direct computation shows that the Casimir operator
\begin{gather}
\label{CasimirWpm}
K_h=\left(\frac{1}{h}W\right)^2-\frac{1}{2}\left(\frac{1}{h}W_h^+ \frac{1}{h}W_h^-+\frac{1}{h}W_h^- \frac{1}
{h}W_h^+\right)
\end{gather}
determines all irreducible unitary representations $\pi_\lambda$ of ${\rm SU}(1,1)$ on the enveloping
algebra $U(\mathfrak{su}(1,1))$ through its eigenvalues $\lambda$.
Moreover, the positive series representation $\pi^+_\lambda$ labelled by $\lambda$ is thus determined by
the set of ladder operators
\begin{gather*}
\pi^+_\lambda\left(\frac{1}{h}W_h^-\right) = E_h^+-E_h^-,
\qquad
\pi^+_\lambda\left(\frac{1}{h}W_h^+\right) = \frac{1}{h}W_h^+,
\qquad
\pi^+_\lambda\left(\frac{1}{h}W\right) = E_h^++\frac{n}{2}I,
\\
\pi^+_\lambda (K_h ) = \left(E_h^++\frac{n}{2}I\right)\left(E_h^++\left(\frac{n}{2}
-1\right)I\right)-\frac{W_h^+}{h}\big(E_h^+-E_h^-\big)
\end{gather*}
underlies the representation space $\ell^2({\mathbb N}_0)$ carrying the family of subspaces
$(\mathcal{H}_{s;h})_{s\in {\mathbb N}_0}$ of the form
\begin{gather*}
\mathcal{H}_{s;h}=\left\{{\bf w}_{s}(x;h)=\left(\frac{1}{h}W_h^+\right)^s{\bf m}_0(x;h)\in\mathcal{P}
 : E_h^+{\bf m}_0(x;h)=E_h^-{\bf m}_0(x;h)=0\right\},
\end{gather*}
while the negative ones is thus determined by the set of ladder operators
\begin{gather*}
\pi^-_\lambda\left(\frac{1}{h}W_h^-\right)=\frac{1}{h}W_h^-,
\qquad
\pi^-_\lambda\left(\frac{1}{h}W_h^+\right)=E_h^+-E_h^-,
\qquad
\pi^-_\lambda\left(\frac{1}{h}W\right)=-E_h^--\frac{n}{2}I,
\\
\pi^-_\lambda (K_h )=\left(E_h^-+\frac{n}{2}I\right)\left(E_h^-+\left(\frac{n}{2}
-1\right)I\right)-\frac{W_h^-}{h}(E_h^+-E_h^-)
\end{gather*}
and underlies the representation space $\ell^2({\mathbb N}_0)$ carrying the family of subspaces
$(\mathcal{H}_{s;-h})_{s\in {\mathbb N}_0}$ of the form
\begin{gather*}
\mathcal{H}_{s;-h}=\left\{{\bf w}_{s}(x;-h)=\left(\frac{1}{h}W_h^-\right)^s{\bf m}_0(x;h)\in\mathcal{P}
 : E_h^+{\bf m}_0(x;h)=E_h^-{\bf m}_0(x;h)=0\right\}.
\end{gather*}

A short computation based on the ladder operator properties~\eqref{su11positive} and~\eqref{su11negative}
shows that the positive/negative discrete series representations $\pi^+_\lambda/\pi^-_\lambda$ of the Lie
group ${\rm SU}(1,1)$ are labeled by the constants $\lambda=\frac{n^2}{4}-\frac{n}{2}-2s$, with $s\in
{\mathbb N}_0$.

\begin{remark}\sloppy
In the context of quantum mechanical systems, this framework may be derived as consequence of a~more
general result~-- the so-called Crum's theorem (cf.~\cite[Subsec\-tions~2.2,~2.3]{OR11}).
\end{remark}

\subsection[The action of ${\rm SO}(n)$ on the lattice]{The action of $\boldsymbol{{\rm SO}(n)}$ on the
lattice}

With the aim of understanding the action of ${\rm SO}(n)$ on the lattice $h{\mathbb Z}^n$ as
a~representation theory carrying canonical generators which are invariant with respect to the orthogonal
Lie algebra $\mathfrak{so}(n)$, we will f\/irst recall some basic concepts and observations about ${\rm
SO}(n)$ and $\mathfrak{so}(n)$ in the context of Clif\/ford algebras.

Let us denote by $\mathcal{B}(x,y)=-\frac{1}{2}(xy+yx)$ the bilinear form generated by the Clif\/ford
vector representations $x=\sum\limits_{j=1}^n x_j {\bf e}_j$ and $y=\sum\limits_{j=1}^n y_j {\bf e}_j$ of
${\mathbb R}^n$.
The set of matrices $T:{\mathbb R}^n \rightarrow {\mathbb R}^n$ with determinant equals 1 for which
$\mathcal{B}(Tx,Ty)=\mathcal{B}(x,y)$ forms a~group under the operation of composition.
This group is called the {\it special orthogonal group of rotations} and it is denoted by~${\rm SO}(n)$.

Let us take a~close look for the ${\rm SO}(n)$-action on the lattice $h{\mathbb Z}^n$ given by the left
regular representation
\begin{gather}
\label{SOnRepresentation}
\Lambda(T){\bf f}(x)={\bf f}\big(T^{-1}x\big)
\qquad
\text{with}
\quad T\in{\rm SO}(n),
\quad
{\bf f}\in\mathcal{P}
\quad
\mbox{and}
\quad
x\in h{\mathbb Z}^n.
\end{gather}

In order to proceed, one can select from ${\rm SO}(n)$ the $1$-parameter subgroups elements $T_{jk}^{\pm
h}(\theta)=\exp\big(\theta S_{jk}^{\pm h}\big)$ generated by exponentiation from the Lie algebra
elements $S_{jk}^{\pm h} \in \mathfrak{so}(n)$ $(1\leq j<k\leq n)$ in such way that each $S_{jk}^{\pm h}$
is skew-symmetric, i.e.\
$S_{jk}^{\pm h}=-S_{kj}^{\pm h}$.

In particular, the canonical elements of the form
\begin{gather}
S_{jk}^{+h}=\mu^{-1}w\left(x_j+\frac{h}{2}\right)\partial_h^{+k}-\mu^{-1}w\left(x_k+\frac{h}{2}
\right)\partial_h^{+j},\nonumber
\\
S_{jk}^{-h}=\mu^{-1}w\left(x_j-\frac{h}{2}\right)\partial_h^{-k}-\mu^{-1}w\left(x_k-\frac{h}{2}
\right)\partial_h^{-j}\label{sonLieAlgebra}
\end{gather}
that correspond to the discrete counterparts of the classical angular momentum operators
$L_{jk}=x_j\partial_{x_k}-x_k\partial_{x_j}$, endow the left representations of ${\rm SO}(n)$ acting on
each subspace $\mathcal{H}_{s;h}$, resp.\ $\mathcal{H}_{s;-h}$ of $\mathcal{P}$, that is
$\Lambda\big(T_{jk}^{+h}(\theta)\big)=\exp\big(\theta S_{jk}^{+h}\big)$, resp.\
$\Lambda\big(T_{jk}^{-h}(\theta)\big)=\exp\big(\theta S_{jk}^{-h}\big)$.

It is straightforwardly to see from~\eqref{Ehpm} that the operator $E_h^+-E_h^-$ commutes with the
skew-symmetric elements~\eqref{sonLieAlgebra} belonging in this way to the center of the enveloping algebra
$U(\mathfrak{so}(n))$.
Then $E_h^+-E_h^-$ commutes with all left regular representations $\Lambda\big(T_{jk}^{+h}(\theta)\big)$,
resp.\ $\Lambda\big(T_{jk}^{-h}(\theta)\big)$, of~${\rm SO}(n)$ and so $E_h^+-E_h^-$ commutes with all
left regular representations of the form~\eqref{SOnRepresentation}.

By virtue of the above Lie algebraic representation we have shown that the family of subspaces
$\left(\mathcal{H}_{s;\pm h}\right)_{s \in {\mathbb N}_0}$ of $\mathcal{P}$ are ${\rm SO}(n)$-invariant on
$h{\mathbb Z}^n$ with respect to $E_h^+-E_h^-$, and so, to all the operators encoded by left-regular
representations of ${\rm SO}(n)$ on $h{\mathbb Z}^n$.

\begin{remark}
The left regular representations $\Lambda\big(T_{jk}^{\pm h}(\theta)\big)$ of ${\rm SO}(n)$ on the
lattice $h{\mathbb Z}^n$ are canonically isomorphic to standard left regular representations
$\Lambda (T_{jk}(\theta) )$ of ${\rm SO}(n)$ on ${\mathbb R}^n$ given in terms of rotations on the
2-dimensional plane with coordinates $(x_j,x_k)$ (cf.~\cite[Chapter~9.1]{VilenkinKlimyk91}).

Indeed, the Clif\/ford-vector-valued extension of the classical Shef\/fer map $\Psi_x$ from ${\mathbb
R}[x]$ to $\mathcal{P}$ def\/ined by linearity from the mapping
\begin{gather*}
 \Psi_x: \ \prod\limits_{j=1}^n x_j^{\alpha_j} \mapsto \prod\limits_{j=1}^n
\big(W_h^{+j}\big)^{\alpha_j}1, \qquad \text{resp.}\qquad  \Psi_x: \ \prod\limits_{j=1}^n x_j^{\alpha_j} \mapsto
\prod\limits_{j=1}^n \big(W_h^{-j}\big)^{\alpha_j}1,
\end{gather*}
satisf\/ies the intertwining relations $\Psi_x x_j=W_{h}^{+j}\Psi_x$ and $\Psi_x
\partial_{x_j}=\partial_{h}^{-j}\Psi_x$, resp.\ $\Psi_x x_j=W_{h}^{-j}\Psi_x$ and $\Psi_x
\partial_{x_j}=\partial_{h}^{+j}\Psi_x$.

This leads to the intertwining relations $\Psi_x L_{jk}=S_{jk}^{+h}\Psi_x$, resp.\ $\Psi_x
L_{jk}=S_{jk}^{-h}\Psi_x$, with $L_{jk}=x_j\partial_{x_k}-x_k\partial_{x_j}$, and hence to the intertwining
property below at the level of ${\rm SO}(n)$:
\begin{gather*}
\Psi_x\Lambda\left(T_{jk}(\theta)\right)=\Lambda\big(T_{jk}^{\pm h}(\theta)\big)\Psi_x.
\end{gather*}

This in turn shows that the 1-parameter representation $T_{jk}^{+h}(\theta)$, resp.\ $T_{jk}^{-h}(\theta)$
of ${\rm SO}(n)$, on~$h{\mathbb Z}^n$ is canonically isomorphic to to the 1-parameter representation
$T_{jk}(\theta)$ of ${\rm SO}(n)$ on~${\mathbb R}^n$.
\end{remark}

\subsection[The Howe dual pair $({\rm SO}(n),\mathfrak{su}(1,1))$]{The Howe dual pair $\boldsymbol{({\rm
SO}(n),\mathfrak{su}(1,1))}$}

Along this section we will study subspaces of Clif\/ford-vector-valued polynomials which are invariant
under the action of the ${\rm SO}(n)\times \mathfrak{su}(1,1)$-module.
The ladder properties~\eqref{su11positive} and~\eqref{su11negative} reveal that for each $s\in {\mathbb
N}_0$ the representation $\pi_\lambda^{+}\left(\frac{1}{h}W_h^+\right)=\frac{1}{h}W_h^+$, resp.
$\pi_\lambda^{-}\left(\frac{1}{h}W_h^-\right)=\frac{1}{h}W_h^-$, maps $\mathcal{H}_{s;h}$, resp.\
$\mathcal{H}_{s;-h}$, into $\mathcal{H}_{s+1;h}$, resp.\ $\mathcal{H}_{s+1;-h}$, while
$\pi_\lambda^{+}\left(\frac{1}{h}W_h^-\right)=\pi_\lambda^{-}\left(\frac{1}{h}W_h^+\right)=E_h^+-E_h^-$
maps $\mathcal{H}_{s;h} \cap \mathcal{H}_{s;-h}$ onto the trivial space $\{ 0\}$.
Also, the positive/negative representations $\pi_\lambda^{\pm}\left(\frac{1}{h}W\right)$, resp.\
$\pi_\lambda^{\pm}(K_h)$, leave the the subspace $\mathcal{H}_{s;h}$, resp.\ $\mathcal{H}_{s;-h}$, invariant.
Here $W$ and $K_h$ denote the multiplication and the Casimir operator labeled by~\eqref{Whpm}
and~\eqref{CasimirWpm}, respectively.

In addition, for any $r=0,1,\ldots,s-1$, the Clif\/ford-vector-valued polynomial spaces
\begin{gather*}
\left(\frac{1}{h}W_h^{\pm}\right)^r\!\left(\mathcal{H}_{s-r;h}\cap\mathcal{H}_{s-r;-h}
\right)=\left\{\left(\frac{1}{h}W_h^{\pm}\right)^r\! {\bf m}_{s-r}(x;h) : {\bf m}_{s-r}(x;h)\in\mathcal{H}
_{s-r;h}\cap\mathcal{H}_{s-r;-h}\right\}
\end{gather*}
are also invariant under the action of $\pi_\lambda^{+}\left(\frac{1}{h}W\right)$ and
$\pi_\lambda^{+}(K_h)$, resp.\
$\pi_\lambda^{-}\left(\frac{1}{h}W\right)$ and $\pi_\lambda^{-}(K_h)$.

So, each ${\rm SO}(n)$-invariant subspace $\mathcal{H}_{s;h}$, resp.\
$\mathcal{H}_{s;-h}$, of $\mathcal{P}$ is isomorphic to the family of subspaces $\left(
\frac{1}{h}W_h^+\right)^r \left(\mathcal{H}_{s-r;h}\cap \mathcal{H}_{s-r;-h}\right)$, resp.\
$\left( \frac{1}{h}W_h^-\right)^r \left(\mathcal{H}_{s-r;h}\cap \mathcal{H}_{s-r;-h}\right)$, with
$r=0,\ldots,s$.
This means that each $\mathcal{H}_{s;h}$, resp.\ $\mathcal{H}_{s;-h}$, appears with inf\/inite multiplicity.

In addition, since $\left( \frac{1}{h}W_h^+\right)^r \left(\mathcal{H}_{s-r;h}\cap
\mathcal{H}_{s-r;-h}\right)=\{ 0\}=\left( \frac{1}{h}W_h^-\right)^r \left(\mathcal{H}_{s-r;h}\cap
\mathcal{H}_{s-r;-h}\right)$ if and only if $s=r$, by virtue of ladder operator
relations~\eqref{su11positive}, resp.~\eqref{su11negative}, one can infer that the direct sum decomposition
\begin{gather*}
\left(\frac{1}{h}W_h^+\right)^r\left(\mathcal{H}_{s-r;h}\cap\mathcal{H}_{s-r;-h}\right)=\mathcal{V}
_h\oplus\mathcal{W}_h,
\\
\text{resp.}
\qquad
\left(\frac{1}{h}W_h^-\right)^r\left(\mathcal{H}_{s-r;h}\cap\mathcal{H}
_{s-r;-h}\right)=\mathcal{V}_h\oplus\mathcal{W}_h,
\end{gather*}
only fulf\/ils for the subspaces $\mathcal{V}_h=\left( \frac{1}{h}W_h^+\right)^r
\left(\mathcal{H}_{s-r;h}\cap \mathcal{H}_{s-r;-h}\right)$ and $\mathcal{W}_h=\{ 0\}$, resp.\
$\mathcal{V}_h=\left( \frac{1}{h}W_h^-\right)^r \left(\mathcal{H}_{s-r;h}\cap
\mathcal{H}_{s-r;-h}\right)$ and $\mathcal{W}_h=\{ 0\}$.
This means that the ${\rm SO}(n)$-invariant subspaces $\left( \frac{1}{h}W_h^+\right)^r
\left(\mathcal{H}_{s-r;h}\cap \mathcal{H}_{s-r;-h}\right)$, resp.\ $\left( \frac{1}{h}W_h^-\right)^r
\left(\mathcal{H}_{s-r;h}\cap \mathcal{H}_{s-r;-h}\right)$, of $\mathcal{P}$ are also irreducible.

In order to collect the inf\/inite multiplicities of $\pi^+_\lambda/\pi^-_\lambda$ carrying the eigenvalues
$\lambda=\frac{n^2}{4}-\frac{n}{2}-2s$ of the Casimir operator~\eqref{CasimirWpm}, it remains to
investigate the $\mathfrak{su}(1,1)$-action on $\mathcal{P}$ regarded as a~${\rm SO}(n)\times
\mathfrak{su}(1,1)$-module that yields the Howe dual pair $({\rm SO}(n),\mathfrak{su}(1,1))$.
For a~sake of readability, the Howe dual pair construction will be only sketched.
Further details arising this construction can be found in~\cite[Chapters 4 \& 5]{GoodmanWallach1998}.

First, notice that the above set of properties produces the following inf\/inite triangle as a~chain
diagram carrying the families of subspaces $\left(\mathcal{H}_{s;h}\right)_{s\in {\mathbb N}_0}$, resp.\
$\left(\mathcal{H}_{s;-h}\right)_{s\in {\mathbb N}_0}$:
\begin{gather*}
\begin{array}{@{}c@{}ccccccccccc@{}}\{0\}&&\mathcal{H}_{0;\pm h}&\leftarrow&\mathcal{H}_{1;\pm h}&\leftarrow&\mathcal{H}_{2;\pm h}
&\leftarrow&\ldots
\\
\;
\\
\{0\}&&\mathcal{H}_{0;h}\cap\mathcal{H}_{0;-h}&\leftarrow&\frac{1}{h}W_h^+\left(\mathcal{H}_{0;h}
\cap\mathcal{H}_{0;-h}\right)&\leftarrow&(\frac{1}{h}W_h^{\pm})^2\left(\mathcal{H}_{0;h}\cap\mathcal{H}
_{0;-h}\right)&\leftarrow&\ldots
\\
&&\oplus&&\oplus&&\oplus&&
\\
&&\{0\}&&\mathcal{H}_{1;h}\cap\mathcal{H}_{1;-h}&\leftarrow&\frac{1}{h}W_h^+\left(\mathcal{H}_{1;h}
\cap\mathcal{H}_{1;-h}\right)&\leftarrow&\ldots
\\
&&&&\oplus&&\oplus&&&&
\\
&&&&\{0\}&&\mathcal{H}_{2;h}\cap\mathcal{H}_{2;-h}&\leftarrow&\ldots
\\
&&&&&&\oplus&&&&
\\
&&&&&&\{0\}&\leftarrow&\ldots
\\
&&&&&&&&&&
\\
&&&&&&&&\ldots
\\
\end{array}
\end{gather*}

In the above triangle diagram, the representations
$\pi_\lambda^+\left(\frac{1}{h}W_h^-\right)=\pi_\lambda^-\left(\frac{1}{h}W_h^+\right)=E_h^+-E_h^-$~--~the
Fourier duals of $\pi_\lambda^\pm \left(\frac{1}{h}W_h^\pm\right)$~--~act as isomorphisms that shift each
individual summand from the right to the left.
The f\/irst line gives the direct sum decomposition of $\mathcal{P}$ in terms of the ${\rm
SO}(n)$-invariant pieces $\mathcal{H}_{s;h}$, resp.
$\mathcal{H}_{s;-h}$, through $h{\mathbb Z}^n$, i.e.
\begin{gather*}
\mathcal{P}=\bigoplus_{s=0}^\infty\mathcal{H}_{s;h}=\bigoplus_{s=0}^\infty\mathcal{H}_{s;-h}.
\end{gather*}

Also, for each $s\in{\mathbb N}_0$ the $(s+1)$-row (which is inf\/inite-dimensional) give rise to
$\mathfrak{su}(1,1)$-modules isomorphic to the ${\rm SO}(n)$-module $\mathcal{H}_{s;h}$, resp.
$\mathcal{H}_{s;-h}$, while the $(s+1)$-column provides the splitting of the subspace $\mathcal{H}_{s;h}$,
resp.\
$\mathcal{H}_{s;-h}$, as a~direct sum in terms of the irreducible pieces $\left( \frac{1}{h}W_h^+\right)^r
\left(\mathcal{H}_{s-r;h}\cap \mathcal{H}_{s-r;-h}\right)$, resp.\ $\left( \frac{1}{h}W_h^-\right)^r
\left(\mathcal{H}_{s-r;h}\cap \mathcal{H}_{s-r;-h}\right)$:
\begin{gather*}
\mathcal{H}_{s;h}=\bigoplus_{r=0}^s\left(\frac{1}{h}W_h^+\right)^r\left(\mathcal{H}_{s-r;h}
\cap\mathcal{H}_{s-r;-h}\right),
\\
\mathcal{H}_{s;-h}=\bigoplus_{r=0}^s\left(\frac{1}{h}W_h^-\right)^r\left(\mathcal{H}_{s-r;h}
\cap\mathcal{H}_{s-r;-h}\right).
\end{gather*}

These chain decompositions lead, in a~multiplicity free way, to the following Fourier decompositions of
$\mathcal{P}$:
\begin{gather}
\mathcal{P}=\bigoplus_{s=0}^\infty\bigoplus_{r=0}^s\left(\frac{1}{h}W_h^+\right)^r\left(\mathcal{H}
_{s-r;h}\cap\mathcal{H}_{s-r;-h}\right),\nonumber\\
\mathcal{P} = \bigoplus_{s=0}^\infty\bigoplus_{r=0}^s\left(\frac{1}{h}W_h^-\right)^r\left(\mathcal{H}
_{s-r;h}\cap\mathcal{H}_{s-r;-h}\right).\label{FourierDecomposition}
\end{gather}
\begin{remark}
The algebra of endomorphisms $\text{End}(\mathcal{P})$ encoded on the Fourier
decompositions~\eqref{FourierDecomposition} is the so-called Weyl algebra of polynomial operators
(cf.~\cite[Chapter~4]{GoodmanWallach1998}).
\end{remark}

\section{Families of special functions}
\label{FamiliesSpecialFunctions}

\subsection{Hypergeometric series representations}

\begin{proposition}
\label{AlmansiDecomposition}
Any sequence of polynomials $\{{\bf m}_s(x;h) : s\in {\mathbb N}_0\}$ satisfying the set of eigenvalue
equations $E_h^+{\bf m}_s(x;h)=E_h^-{\bf m}_s(x;h)=s{\bf m}_s(x;h)$ is determined in a~unique way as
\begin{gather*}
{\bf m}_s(x;h)=\gamma_{s}{\bf w}_{s}(x;\pm h)
\qquad
\text{with}
\quad
{\bf w}_{s}(x;\pm h)\in\mathcal{H}_{s;\pm h}
\quad
\text{and}
\quad
\gamma_s\in{\mathbb R}.
\end{gather*}
Moreover, the constants $\gamma_s$ are given by
\begin{gather*}
\gamma_s=\sum_{r=0}^s\frac{(-1)^r(-s-n-1)_r}{(-2s-n+2)_{r}}\left(
\begin{matrix}s
\\
r
\end{matrix}
\right).
\end{gather*}
Hereby $(a)_r=a(a+1)\cdots (a+r-1)$ denotes the Pochhammer symbol.
\end{proposition}

\begin{proof} To prove this, let us consider the family $\left\{ {\bf m}_s(x;h) : s\in {\mathbb N}_0\right\}$ of
Clif\/ford-vector-valued polynomials, each of them given as a~solution of the eigenvalue equations
\begin{gather}
\label{EigenvalueConstraint}
E_h^+{\bf m}_s(x;h)=E_h^-{\bf m}_s(x;h)=s{\bf m}_s(x;h).
\end{gather}

From the Fourier decompositions labeled by~\eqref{FourierDecomposition} each ${\bf w}_s(x;\pm h)\in
\mathcal{H}_{s;h}$ may be written in a~unique way as
\begin{gather*}
{\bf w}_s(x;\pm h)=\sum_{r=0}^s\left(\frac{1}{h}W_h^\pm\right)^r{\bf m}_{s-r}(x;h).
\end{gather*}

From the ladder operator relations~\eqref{su11positive} and~\eqref{su11negative}, one obtain, for any
$r=0,1,\ldots,s$, the mapping property
$\left(\frac{1}{h}W_h^\pm\right)^r\left(E_h^+-E_h^-\right)^r:\mathcal{H}_{s;\pm h} \rightarrow
\mathcal{H}_{s;\pm h}$.
Then, for a~given polynomial sequence $\left\{ {\bf w}_s(x;\tau) : s\in {\mathbb N}_0\right\}$ of
$\mathcal{P}$, one can compute each ${\bf m}_s(x;h)$ from the formulae
\begin{gather}
\label{msPolynomials}
{\bf m}_s(x;h)=\sum_{r=0}^s c_{r,s}\left(\frac{1}{h}W_h^\pm\right)^r\widetilde{{\bf w}}_{s-r}(x;\pm h).
\end{gather}
where $\widetilde{{\bf w}}_{s-r}(x;\pm h)=\left(E_h^+-E_h^-\right)^r {\bf w}_{s}(x;\pm h)$ and the
constants $c_{r,s}\in {\mathbb R}$ are determined from the constraint~\eqref{EigenvalueConstraint}.

Combination of the graded commuting relation provided from Proposition~\ref{su11EhpmProposition}:
\begin{gather*}
\left[E_h^+-E_h^-,\left(\frac{1}{h}W_h^\pm\right)^r\right]=r(2E_h^\pm+(n-r+1)I)\left(\frac{1}{h}
W_h^\pm\right)^{r-1}
\end{gather*}
with the eigenvalue properties
\begin{gather*}
E_h^\pm\left[\left(\frac{1}{h}W_h^\pm\right)^{r-1}\widetilde{{\bf w}}_{s-r}
(x;\pm h)\right]=(s-1)\left(\frac{1}{h}W_h^\pm\right)^{r-1}\widetilde{{\bf w}}_{s-r}(x;\pm h)
\end{gather*}
lead to the recursive relations
\begin{gather*}
\big(E_h^+-E_h^-\big)\left[\left(\frac{1}{h}W_h^\pm\right)^r\widetilde{{\bf w}}_{s-r}(x;\pm h)\right]
\\
\qquad {} =r(2s+n-r-1)\left(\frac{1}{h}W_h^\pm\right)^{r-1}\widetilde{{\bf w}}_{s-r}(x;\pm h)+\left(\frac{1}{h}
W_h^\pm\right)^r\widetilde{{\bf w}}_{s-r-1}(x;\pm h)
\end{gather*}
and moreover, to the following linear expansions
\begin{gather*}
E_h^+{\bf m}_{s}(x;h)-E_h^-{\bf m}_{s}(x;h)=\sum_{r=1}^s d_{r,s}\left(\frac{1}{h}
W_h^\pm\right)^r\widetilde{{\bf w}}_{s-r-1}(x;\pm h),
\end{gather*}
where the coef\/f\/icients $d_{r,s}$ are given by $d_{r,s}=(r+1)(2s+n-r-2)c_{r+1,s}+c_{r,s}$.

Hence ${\bf m}_{s}(x;h)$ satisf\/ies the constraint~\eqref{EigenvalueConstraint} if and only if $d_{r,s}=0$
holds for every $r=0,1,\ldots,s$, that is, each $c_{r,s}$ is determined from the condition $c_{0,s}=1$ and
from the constraint
\begin{gather*}
c_{r+1,s}=\frac{c_{r,s}}{(r+1)(-2s-n+r+2)}.
\end{gather*}
Therefore
\begin{gather*}
c_{r,s}=\prod_{q=1}^{r}\frac{1}{q(-2s-n+q+1)}=\frac{1}{r!(-2s-n+2)_r}.
\end{gather*}

Now it remains to show the relations ${\bf m}_s(x;h)=\gamma_s {\bf w}_{s}(x;\pm h)$.
Iterating $r$ times the operator $E_h^+-E_h^-$, one obtain from~\eqref{su11positive} the recursive
relations
\begin{gather*}
\left(E_h^+-E_h^-\right)^r{\bf w}_{s}(x;h) = (-1)^r\frac{s!}{(s-r)!}(-s-n)_r{\bf w}_{s-r}(x;h),
\\
\left(\frac{1}{h}W_h^+\right)^r{\bf w}_{s-r}(x;h) = {\bf w}_{s}(x;h).
\end{gather*}

Then we have $\left(\frac{1}{h}W_h^+\right)^r\widetilde{{\bf w}}_{s-r}(x;\pm h)=(-1)^r\frac{s!}{(s-r)!}
(-s-n-1)_r{\bf w}_{s}(x;\pm h)$, and therefore, the set of relations~\eqref{msPolynomials} are equivalent
to ${\bf m}_s(x;h)=\gamma_s {\bf w}_{s}(x;\pm h)$, with
\begin{gather*}
\gamma_s=\sum_{r=0}^s\frac{(-1)^r(-s-n-1)_r}{(-2s-n+2)_{r}}\left(
\begin{matrix}s
\\
r
\end{matrix}
\right).\tag*{\qed}
\end{gather*}
\renewcommand{\qed}{}
\end{proof}

\begin{remark}
A simple computation involving the binomial identity $(-1)^r\left(
\begin{matrix} s
\\
r
\end{matrix}
\right)=\frac{(-s)_r}{r!}$ shows that the constant $\gamma_{s}$ corresponds to the $s$-term truncation of
the hypergeometric function ${}_2F_1(a,b;c;z)=\sum\limits_{r=0}^\infty \frac{(a)_r (b)_r}{(c)_r}\frac{z^r}{
r!}$ labelled by the parameters $a=-s-n-1$, $b=-s$, $c=-2s-n+2$ and $z=1$.
\end{remark}

\subsection{Application to Cauchy problems}

On this section it will be studied families of Clif\/ford-vector-valued polynomials given as solutions of
the following homogeneous Cauchy problem in $[0,\infty)\times h{\mathbb Z}^n$:
\begin{alignat}{3}
& \partial_t{\bf g}(t,x)+E_h^+{\bf g}(t,x)-E_h^-{\bf g}(t,x)=0,\qquad && t>0, &\nonumber\\
& {\bf g}(0,x)={\bf f}(x),\qquad && t=0,&\label{CauchyProblem}\\
& E_h^+{\bf g}(t,x)=E_h^-{\bf g}(t,x),\qquad && t\geq0. &\nonumber
\end{alignat}

The solution of the Cauchy problem~\eqref{CauchyProblem} may be written formally as ${\bf
g}(t,x)=\mathbb{E}_h(t){\bf f}(x)$ with ${\bf f}(x)\in \bigoplus_{s=0}^\infty\mathcal{H}_{s;h} \cap
\mathcal{H}_{s;-h}$.
Hereby $\mathbb{E}_h(t)=\exp(tE_h^--tE_h^+)$ is a~one-parameter representation of the Lie group ${\rm
SU}(1,1)$.

Note that the graded commuting property $\left[E_h^+-E_h^-,\mathbb{E}_h(t)\right]=0$ for each $t \geq 0$
assures that $\left( \mathbb{E}_h(t)\right)_{t\geq 0}$ is a~semigroup, i.e.
\begin{gather*}
\mathbb{E}_h(0)=I
\qquad
\text{and}
\qquad
\mathbb{E}_h(t+\tau)=\mathbb{E}_h(t)\mathbb{E}_h(\tau) \quad \mbox{for all} \ \ t,\tau\geq0.
\end{gather*}

The next sequence of results will makes clear that, for any $t\geq 0$, the operator $\mathbb{E}_h(t)$
leaves the space of Clif\/ford-vector-valued polynomials $\mathcal{P}$ invariant.
\begin{lemma}
\label{EhpmWpmIntertwiningDiscreteSeries}
We have the following intertwining properties carrying the semigroup operator
$\mathbb{E}_h(t)=\exp(tE_h^--tE_h^+)$:
\begin{gather*}
\big(tE_h^-+(1-t)E_h^+\big)\mathbb{E}_h(t)=\mathbb{E}_h(t)E_h^+,
\\
\left(\frac{1}{h}W_h^+-t\left(E_h^++E_h^-+nI\right)\right)\mathbb{E}_h(t)=\mathbb{E}_h(t) \frac{1}{h}W_h^+.
\end{gather*}
\end{lemma}

\begin{proof} In order to start proving the intertwining property
\begin{gather*}
\big(tE_h^-+(1-t)E_h^+\big)\mathbb{E}_h(t)=\mathbb{E}_h(t)E_h^+,
\end{gather*}
one can start to compute, for each $s \in {\mathbb N}$, the graded commutator $\left[2E_h^+
+nI,\left(E_h^--E_h^+\right)^s\right]$.
A~short computation based on Lemma~\ref{su11EhpmLemma} shows that
\begin{gather*}
\left[E_h^+,E_h^--E_h^+\right]=\left[E_h^++\frac{n}{2}I,E_h^--E_h^+\right]=E_h^+-E_h^-.
\end{gather*}

Thus, direct application of Lemma~\ref{ABsLemma} for $A=E_h^+$ and $B=E_h^--E_h^+$ results into
\begin{gather*}
\left[E_h^+,\left(E_h^--E_h^+\right)^s\right]=s\left(E_h^+-E_h^-\right)\left(E_h^--E_h^+\right)^{s-1}.
\end{gather*}

This leads to
\begin{gather*}
\left[E_h^+,\mathbb{E}_h(t)\right]=\sum_{r=0}^\infty\frac{t^r}{r!}
\left[E_h^+,\left(E_h^--E_h^+\right)^r\right]
\\
\hphantom{\left[E_h^+,\mathbb{E}_h(t)\right]}{}
=\sum_{r=1}^\infty\frac{t^r}{(r-1)!}(E_h^+-E_h^-)\left(E_h^--E_h^+\right)^{r-1}
=\left(tE_h^+-tE_h^-\right)\mathbb{E}_h(t).
\end{gather*}

Therefore, the above graded commuting relation is equivalent to the intertwining property
\begin{gather*}
\mathbb{E}_h(t)E_h^+=(tE_h^-+(1-t)E_h^+)\mathbb{E}_h(t).
\end{gather*}

For the proof of $\left(\frac{1}{h}W_h^+-t\left(E_h^+ +E_h^- +
nI\right)\right)\mathbb{E}_h(t)=\mathbb{E}_h(t) \frac{1}{h}W_h^+$, one needs to compute, for every $r \in
{\mathbb N}$, the graded commutator $\left[\frac{1}{h}W_h^+,\left(E_h^+-E_h^-\right)^r\right]$ based on
Lemma~\ref{ABsLemma}.
First, recall that the relation $[E_h^--E_h^+,E_h^++\frac{n}{2}I]=E_h^--E_h^+$ that follows from
Lemma~\ref{su11EhpmLemma} is equivalent to the intertwining property
\begin{gather*}
\left(E_h^--E_h^+\right)\left(E_h^++\frac{n}{2}I\right)=\left(E_h^++\left(\frac{n}{2}
+1\right)I\right)\left(E_h^--E_h^+\right).
\end{gather*}

Induction over $r\in {\mathbb N}$ gives
\begin{gather*}
\left(E_h^--E_h^+\right)^r\left(E_h^++\frac{n}{2}I\right)=\left(E_h^++\left(\frac{n}{2}
+r\right)I\right)\left(E_h^--E_h^+\right)^r.
\end{gather*}

On the other hand, from $\left[ \frac{1}{h}W_h^+,E_h^--E_h^+\right]=2E_h^++nI$ (see
Lemma~\ref{su11EhpmLemma}) and from direct application of Lemma~\ref{ABsLemma} results into the following
set of relations carrying the generators $A=\frac{1}{h}W_h^+$ and $B=E_h^--E_h^+$:
\begin{gather*}
\left[\frac{1}{h}W_h^+,\left(E_h^--E_h^+\right)^s\right]=\sum_{r=0}^{s-1}
\left(E_h^--E_h^+\right)^r \left(2E_h^++nI\right) \left(E_h^--E_h^+\right)^{s-1-r}
\\
\hphantom{\left[\frac{1}{h}W_h^+,\left(E_h^--E_h^+\right)^s\right]}{}
=\sum_{r=0}^{s-1}\left(2E_h^++\left(n+2r\right)I\right)\left(E_h^--E_h^+\right)^{s-1}\\
\hphantom{\left[\frac{1}{h}W_h^+,\left(E_h^--E_h^+\right)^s\right]}{}
=s\left(2E_h^++(n+s-1)I\right)\left(E_h^+-E_h^-\right)^{s-1}.
\end{gather*}

This leads to
\begin{gather*}
\left[\frac{1}{h}W_h^+,\mathbb{E}_h(t)\right]=\sum_{s=0}^\infty\frac{t^s}{s!}\left[\frac{1}{h}
W_h^+,\left(E_h^--E_h^+\right)^s\right]
=\sum_{s=1}^\infty\frac{t^s}{(s-1)!}(2E_h^++nI)\left(E_h^--E_h^+\right)^{s-1}
\\
\hphantom{\left[\frac{1}{h}W_h^+,\mathbb{E}_h(t)\right]=}
{}+\sum_{s=2}^\infty\frac{t^s}{(s-2)!}\left(E_h^--E_h^+\right)^{s-1}
=t\left(E_h^++E_h^-+nI\right)\mathbb{E}_h(t).
\end{gather*}
and hence, to the intertwining property
\begin{gather*}
\left(\frac{1}{h}W_h^+-t\left(E_h^++E_h^-+nI\right)\right)\mathbb{E}_h(t)=\mathbb{E}_h(t) \frac{1}{h}W_h^+.\tag*{\qed}
\end{gather*}
\renewcommand{\qed}{}
\end{proof}
\begin{proposition}
For any $t\geq 0$ and $s\in {\mathbb N}_0$, we have the mapping property
\begin{gather*}
\mathbb{E}_h(t):\ \mathcal{H}_{s;h}\cap\mathcal{H}_{s;-h}\rightarrow\mathcal{H}_{s;h}\cap\mathcal{H}_{s;-h}.
\end{gather*}
Moreover the semigroup $(\mathbb{E}_h(t))_{t\geq 0}$ leaves invariant the space of Clifford-vector-valued
polynomials $\mathcal{P}$.
\end{proposition}

\begin{proof}
Recall that from the Fourier decomposition~\eqref{FourierDecomposition} any ${\bf f}(x)\in \mathcal{P}$ may
be written uniquely as
\begin{gather*}
{\bf f}(x)=\sum_{s=0}^\infty\sum_{r=0}^s\left(\frac{1}{h}W_h^+\right)^s{\bf m}_{s-r}(x;h),
\end{gather*}
where $\left\{{\bf m}_s(x;h) : s\in {\mathbb N}\right\}$ is a~sequence of polynomials satisfying the set of
eigenvalue equations $ E_h^+ {\bf m}_s(x;h)=E_h^- {\bf m}_s(x;h)=s{\bf m}_s(x;h)$.

From Lemma~\ref{EhpmWpmIntertwiningDiscreteSeries}, one can easily see that the operator $\mathbb{E}_h(t)$
intertwines $E_h^+$ and $tE_h^-+(1-t)E_h^+$.
Since $E_h^+-E_h^-$ commutes with $\mathbb{E}_h(t)$, it follows that the set of functions ${\bf
m}_s(t,x;h):=\mathbb{E}_h(t){\bf m}_s(x;h)$ satisfy the coupled system of equations
\begin{gather*}
\left(tE_h^-+(1-t)E_h^+\right){\bf m}_s(t,x;h) = s{\bf m}_s(t,x;h),
\qquad
\left(E_h^+-E_h^-\right){\bf m}_s(t,x;h) = 0.
\end{gather*}

From the above system of equations, one can infer the eigenvalue relations $ E_h^+ {\bf m}_s(t,x;h)=E_h^-
{\bf m}_s(t,x;h)=s{\bf m}_s(t,x;h)$.
Thus we have shown the mapping property $\mathbb{E}_h(t):\mathcal{H}_{s;h} \cap \mathcal{H}_{s;-h}
\rightarrow \mathcal{H}_{s;h} \cap \mathcal{H}_{s;-h}$.

Finally, since the operator $\frac{1}{h}W_h^+-t\left(E_h^++E_h^-+nI\right)$ is canonically isomorphic to
$\frac{1}{h}W_h^+$, one can see that
\begin{gather*}
\mathbb{E}_h(t) \left[\left(\frac{1}{h}W_h^+\right)^r{\bf m}_{s-r}(x;h)\right]=\left(\frac{1}{h}
W_h^+-t\left(E_h^++E_h^-+nI\right)\right)^r{\bf m}_{s-r}(t,x;h).
\end{gather*}

From the representation of $E_h^++E_h^-+nI$ given by~\eqref{EhpmLinearCombinations}, one can conclude that
the right-hand side of the above relation is a~Clif\/ford-vector-valued polynomial but also an
eigenfunction of $tE_h^-+(1-t)E_h^+$ carrying the eigenvalue $s\in {\mathbb N}_0$.

Therefore $\mathbb{E}_h(t) {\bf f}(x)=\sum\limits_{s=0}^\infty \sum\limits_{r=0}^s \mathbb{E}_h(t)
\left[\left(\frac{1}{h}W_h^+\right)^r{\bf m}_{s-r}(x;h)\right] \in \mathcal{P}$.
\end{proof}

\begin{remark}
One can see from a~direct combination of Lemma~\ref{EhpmWpmIntertwiningDiscreteSeries}
with~\eqref{EhpmLadderWj} that for $t=0$, the polynomial solution provided from the initial condition ${\bf
g}(0,x)={\bf f}(x)$ belongs to $\bigoplus_{s=0}^\infty \mathcal{H}_{s;h}$ while for $t=1$ the solution
${\bf g}(1,x)=\mathbb{E}_h(1){\bf f}(x)$ belongs to $\bigoplus_{s=0}^\infty \mathcal{H}_{s;-h}$.
So, the semigroup $(\mathbb{E}_h(t))_{t\geq 0}$ gives, in particular, a~direct link between the positive
series representation of ${\rm SU}(1,1)$ with the negative ones.
\end{remark}
\begin{remark}
One can see from Lemma~\ref{EhpmWpmIntertwiningDiscreteSeries} that
$\frac{1}{2}\left(\frac{1}{h}W_h^++\frac{1}{h}W_h^-\right)$ is the Fourier dual of $E_h^+-E_h^-$ carrying
the parameter $t=\frac{1}{2}$.

One can also deduce from Lemma~\ref{EhpmWpmIntertwiningDiscreteSeries} that the solutions of the eigenvalue
problem
\begin{gather*}
E_h^+{\bf g}_s (t,x )+E_h^-{\bf g}_s (t,x )=2s{\bf g}_s (t,x )
\end{gather*}
correspond to ${\bf g}_s\left(\frac{1}{2},x\right)=\mathbb{E}_h\left(\frac{1}{2}\right){\bf f}(x)$, with
${\bf f}(x) \in \mathcal{H}_{s;h}$.
Moreover, given ${\bf g}_s(1,x) \in \mathcal{H}_{s;-h}$, one can compute ${\bf
g}_s\left(\frac{1}{2},x\right)$ by letting act the inverse of $\mathbb{E}_h\left(\frac{1}{2}\right)$ on
${\bf g}_s(1,x)$, i.e.
\begin{gather*}
{\bf g}_s\left(\frac{1}{2},x\right)=\mathbb{E}_h\left(-\frac{1}{2}\right){\bf g}_s(1,x).
\end{gather*}
\end{remark}

The next corollary, that follows straightforwardly from the combination of the above proposition with
Proposition~\ref{AlmansiDecomposition}, fully characterize the Clif\/ford-vector-valued polynomial
solutions of the Cauchy problem~\eqref{CauchyProblem} as hypergeometric ${}_0F_1$-series expansions.
\begin{corollary}
Let ${\bf f}(x)\in \bigoplus_{s=0}^\infty \mathcal{H}_{s,h} \cap \mathcal{H}_{s,-h}$.
If there is a~sequence of polynomials $\{{\bf w}_s(x;h):\, s\in {\mathbb N}_0 \}$, each of them belonging
to $\mathcal{H}_{s,h}$, then the solutions ${\bf g}(t,x)=\mathbb{E}_h(t){\bf f}(x)$ of the Cauchy
problem~\eqref{CauchyProblem} are given explicitly in terms of the hypergeometric ${}_0F_1$-series
expansion
\begin{gather*}
{\bf f}(x)=\sum_{s=0}^\infty\left[{}_0F_1(-2s-n+2;\partial_t)t^{s}\right]_{t=1} {\bf w}_s(x;h),
\end{gather*}
where ${}_0F_1(c;z)$ denotes the hypergeometric function ${}_0F_1(c;z)=\sum\limits_{r=0}^\infty
\frac{1}{(c)_r}\frac{z^r}{r!}$.
\end{corollary}

\begin{proof} From Proposition~\ref{AlmansiDecomposition} it follows that each ${\bf f}(x)\in
\bigoplus_{s=0}^\infty \mathcal{H}_{s,h} \cap \mathcal{H}_{s,-h}$ may be written as
\begin{gather*}
{\bf f}(x)=\sum_{s=0}^\infty\gamma_s {\bf w}_s(x;h)
\qquad
\text{with}
\quad
\gamma_s=\sum_{r=0}
^s\frac{(-1)^r(-s-n-1)_r}{(-2s-n+2)_{r}}\left(
\begin{matrix}s
\\
r
\end{matrix}
\right).
\end{gather*}

A short computation based on the lowering properties $(\partial_t)^r t^{s}=\frac{s!}{(s-r)!}t^{s-r}$ for
$s\geq r$ and $(\partial_t)^r t^{s}=0$ for $s<r$ yields $ \gamma_s=\left[{}_0F_1(-2s-n+2;\partial_t)
t^{s}\right]_{t=1}$.
\end{proof}

\subsection*{Acknowledgements}

The work of the author was supported by the fellowship 13/07590-8 of FAPESP (S.P., Brazil) and by the
project PTDC/MAT/114394/2009 funded by FCT (Portugal) through the European program COMPETE/FEDER.
The author would like to thank to the anonymous referees for the careful reading and for the criticism
through the reports.
This allows to improve the quality of the former manuscript in a~clever style and for Waldyr A.~Rodrigues
Jr.\ (IMECC--Unicamp) for the careful reading of the f\/inal version and for the important remarks
concerning Section~\ref{ScopeProblemsSection}.
The major part of this work was developed when the author was a~research member of the Centre for
Mathematics from University of Coimbra (Portugal).
The author wish to express along these lines his gratitude to all members of the research centre for the
excellent atmosphere and for the constant support that made these last three years a~pure enjoyment.

\pdfbookmark[1]{References}{ref}
\LastPageEnding


\begin{thebibliography}{99}
\footnotesize\itemsep=0pt

\bibitem{Cartier00}
Cartier P., Mathemagics (a tribute to {L}.~{E}uler and {R}.~{F}eynman), in
  Noise, Oscillators and Algebraic Randomness ({C}hapelle des {B}ois, 1999),
  \href{http://dx.doi.org/10.1007/3-540-45463-2_2}{\textit{Lecture Notes in Phys.}}, Vol.~550, Springer, Berlin, 2000, 6--67.

\bibitem{RSKS10}
De~Ridder H., De~Schepper H., K{\"a}hler U., Sommen F., Discrete function
  theory based on skew {W}eyl relations, \href{http://dx.doi.org/10.1090/S0002-9939-2010-10480-X}{\textit{Proc. Amer. Math. Soc.}}
  \textbf{138} (2010), 3241--3256.

\bibitem{RSS12}
De~Ridder H., De~Schepper H., Sommen F., Fueter polynomials in discrete
  {C}lif\/ford analysis, \href{http://dx.doi.org/10.1007/s00209-011-0932-5}{\textit{Math.~Z.}} \textbf{272} (2012), 253--268.

\bibitem{DSS92}
Delanghe R., Sommen F., Sou{\v{c}}ek V., Clif\/ford algebra and spinor-valued
  functions. A~function theory for the Dirac operator, \href{http://dx.doi.org/10.1007/978-94-011-2922-0}{\textit{Mathematics and
  its Applications}}, Vol.~53, Kluwer Academic Publishers Group, Dordrecht,
  1992.

\bibitem{BLR98}
Di~Bucchianico A., Loeb D.E., Rota G.C., Umbral calculus in {H}ilbert space, in
  Mathematical Essays in Honor of {G}ian-{C}arlo {R}ota ({C}ambridge, {MA},
  1996), \textit{Progr. Math.}, Vol.~161, Birkh\"auser Boston, Boston, MA,
  1998, 213--238.

\bibitem{DHS96}
Dimakis A., M{\"u}ller-Hoissen F., Striker T., Umbral calculus, discretization,
  and quantum mechanics on a~lattice, \href{http://dx.doi.org/10.1088/0305-4470/29/21/017}{\textit{J.~Phys.~A: Math. Gen.}}
  \textbf{29} (1996), 6861--6876, \href{http://arxiv.org/abs/quant-ph/9509014}{quant-ph/9509014}.

\bibitem{Eelbode12}
Eelbode D., Monogenic {A}ppell sets as representations of the {H}eisenberg
  algebra, \href{http://dx.doi.org/10.1007/s00006-012-0330-z}{\textit{Adv. Appl. Clifford Algebr.}} \textbf{22} (2012), 1009--1023.

\bibitem{FK07}
Faustino N., K{\"a}hler U., Fischer decomposition for dif\/ference {D}irac
  operators, \href{http://dx.doi.org/10.1007/s00006-006-0016-5}{\textit{Adv. Appl. Clifford Algebr.}} \textbf{17} (2007), 37--58,
  \href{http://arxiv.org/abs/math.CV/0609823}{math.CV/0609823}.

\bibitem{FR11}
Faustino N., Ren G., ({D}iscrete) {A}lmansi type decompositions: an umbral
  calculus framework based on {${\mathfrak{osp}}(1|2)$} symmetries,
  \href{http://dx.doi.org/10.1002/mma.1498}{\textit{Math. Methods Appl. Sci.}} \textbf{34} (2011), 1961--1979,
  \href{http://arxiv.org/abs/1102.5434}{arXiv:1102.5434}.

\bibitem{GoodmanWallach1998}
Goodman R., Wallach N.R., Representations and invariants of the classical
  groups, \textit{Encyclopedia of Mathematics and its Applications}, Vol.~68,
  Cambridge University Press, Cambridge, 1998.

\bibitem{Howe89TransAMS}
Howe R., Remarks on classical invariant theory, \href{http://dx.doi.org/10.2307/2001418}{\textit{Trans. Amer. Math.
  Soc.}} \textbf{313} (1989), 539--570.

\bibitem{LTW04}
Levi D., Tempesta P., Winternitz P., Umbral calculus, dif\/ference equations and
  the discrete {S}chr\"odinger equation, \href{http://dx.doi.org/10.1063/1.1780612}{\textit{J.~Math. Phys.}} \textbf{45}
  (2004), 4077--4105, \href{http://arxiv.org/abs/nlin.SI/0305047}{nlin.SI/0305047}.

\bibitem{MT08}
Malonek H.R., Tomaz G., Bernoulli polynomials and {P}ascal matrices in the
  context of {C}lif\/ford analysis, \href{http://dx.doi.org/10.1016/j.dam.2008.06.009}{\textit{Discrete Appl. Math.}} \textbf{157}
  (2009), 838--847.

\bibitem{OR11}
Odake S., Sasaki R., Discrete quantum mechanics, \href{http://dx.doi.org/10.1088/1751-8113/44/35/353001}{\textit{J.~Phys.~A: Math.
  Theor.}} \textbf{44} (2011), 353001, 47~pages, \href{http://arxiv.org/abs/1104.0473}{arXiv:1104.0473}.

\bibitem{Sommen97}
Sommen F., An algebra of abstract vector variables, \textit{Portugal. Math.}
  \textbf{54} (1997), 287--310.

\bibitem{Tempesta08}
Tempesta P., On {A}ppell sequences of polynomials of {B}ernoulli and {E}uler
  type, \href{http://dx.doi.org/10.1016/j.jmaa.2007.07.018}{\textit{J.~Math. Anal. Appl.}} \textbf{341} (2008), 1295--1310.

\bibitem{VilenkinKlimyk91}
Vilenkin N.J., Klimyk A.U., Representation of {L}ie groups and special
  functions. {V}ol.~1.~Simplest Lie groups, special functions and integral
  transforms, \textit{Mathematics and its Applications (Soviet Series)},
  Vol.~72, Kluwer Academic Publishers Group, Dordrecht, 1991.

\end{thebibliography}
\end{document}